\documentclass[12pt]{article}
\usepackage{bezier,latexsym}
\usepackage{epsfig,afterpage}
\usepackage{amsmath, amssymb}

\usepackage{color}
\catcode`\@=11 
\@addtoreset{equation}{section} 
\catcode`\@=12 
\newcommand\real{{\rm I\! R}}
\newcommand\nat{{\rm I\! N}}
\newcommand\beweisende{\ \hfill$\Box$\break}

\newtheorem{theorem}{Theorem}[section]

\newcommand\proof{{\noindent \it Proof.\ }}

\def\bql#1{\begin{equation}\label{#1}}
\def\bq{\begin{equation}}
\def\eq{\end{equation}}
\def\eref#1{{\rm (\ref{#1})}}
\newcommand{\fet}[1]{\mbox{\boldmath ${#1}$}}
\usepackage[colorlinks=true, linkcolor=blue, citecolor =green, pagebackref=false, bookmarks=false, pdfpagemode=UseNone]{hyperref}

\usepackage{array,multirow}
\usepackage{hhline}
{
\usepackage[table]{xcolor}

\begin{document}

\title{On the Non-Linear Integral Equation Approach for an \\ Inverse Boundary Value 
Problem \\ for the Heat Equation}

\author{Roman Chapko\thanks{Faculty of Applied Mathematics and Informatics, Ivan Franko National University of Lviv, 79000
Lviv, Ukraine}\, and Leonidas Mindrinos\thanks{Johann Radon Institute for Computational and Applied Mathematics (RICAM), Austrian Academy of Sciences, Linz, Austria}
}

\date{ }
\maketitle

\begin{abstract}
We consider the inverse problem of reconstructing the interior boundary curve of a doubly connected domain from the knowledge of the temperature and the thermal flux on the exterior boundary curve.  The use of the Laguerre transform in time leads to a sequence of stationary inverse problems. Then, the application of the modified single-layer ansatz, reduces the problem to a sequence of systems of non-linear boundary integral equations. An iterative algorithm is developed for the numerical solution of the obtained integral equations.  We find the Fr\'echet derivative of the corresponding integral operator and we show the unique solvability of the linearized equation. Full discretization is realized by a trigonometric quadrature method. Due to the inherited ill-possedness of the derived system of linear equations we apply the Tikhonov regularization. 
The numerical results show that the proposed method produces accurate and stable reconstructions.

\vspace*{1cm} {\em Keywords:} doubly connected domain; boundary reconstruction; Laguerre transform; modified single layer potentials; non-linear boundary integral equations; tri\-go\-nometric quadrature method; Newton method; Tikhonov regularization.
\end{abstract}

\section{Introduction}

Integral equation methods are among the most frequently used and effective methods for solving problems in applied sciences. The main advantage is that they are methods for dimension reduction and instead of solving a problem in an
unbounded region, we obtain the solution from a boundary integral equation. This has many applications to a wide class of direct boundary value problems but also to various inverse problems. 
 
In this work, we deal with the problem of reconstructing the inner boundary of a medium. There are different ways to solve this problem. Once the inverse problem is written as a non-linear operator equation, then its approximate solution can be found by Newton's method (see i.e. \cite{ChKr3,ChKr2,CoKr}). 
The inverse problem can also be directly reduced to a system of non-linear integral equations. This is the reciprocity gap approach, based on the Green's integral theorem. In \cite{IvJo,IvKr,KrRu}, it was used for the case of the Laplace and Helmholtz equations. 
Another way to reduce the problem of reconstructing the boundary curve of a medium to a system of non-linear integral equations is to use potentials for the integral representation of the solution \cite{AltKre12,Greece,ChIv2,GinMin19,JohSle07}. 
Representing the unknown function in a form of, for example, a single-layer potential, the inverse problem is reduced to a system of non-linear equations with respect to the unknown densities and the unknown boundary of the domain. Then, we may apply an iterative scheme.

In the case of non-linear inverse problem for the non-stationary heat conductivity equation, an approach with the use of the indirect method of integral equations for the solution of the corresponding non-linear operator equation was developed in \cite{ChKr1,ChKr2}.  In this paper, we want to use our experience with solving non-stationary direct problems by a combination of the Laguerre transform with integral equations \cite{ChJo,ChKr,ChKr3}.
The application of the Laguerre transform with respect to the time variable, results to the semi-discretization of the given non-linear inverse parabolic problem. 
This allows us to obtain a sequence of stationary inverse problems, which is then further reduced by the  method of integral potentials to a sequence of systems of non-linear integral equations.

The paper is organized as follows. In the following, we present the governing equations and we formulate the corresponding inverse problem. In \autoref{Sec2}, 
 we apply the Laguerre transform to the inverse boundary problem with respect to the time variable and we reduce it to a sequence of stationary boundary problems. Then, using the modified single-layer ansatz, we rewrite it equivalently as a sequence of non-linear boundary integral equations. Here, the proposed iterative scheme for its solution is also described. The \autoref{sec3} contains the numerical solution of system of the integral equations. Firstly, we consider the well-posed system of the ``field" equations and we apply the Nystr\"om method for its discretization. Then, we linearize the ill-posed ``data" equation and apply the collocation method for its discretization. The Tikhonov regularization is used  for solving the derived ill-conditioned linear system.  Finally, in \autoref{sec4} we provide numerical examples where the reconstruction of the unknown boundary is feasible for different setups.

Let $\Omega$ be a doubly connected domain in $\real^2$ with smooth boundary $\Gamma$ of class $C^2.$
We assume that $\Gamma$ consists of two disjoints curves $\Gamma_1$ and $\Gamma_2$, meaning $\Gamma = \Gamma_1 \cup \Gamma_2,$ with $\Gamma_1 \cap \Gamma_2 = \emptyset,$ such that $\Gamma_1$ is contained in the interior of $\Gamma_2.$

We consider the following initial boundary value problem for the heat equation
\bql{1.1}
\frac{1}{\alpha}\frac{\partial u}{\partial t} =\Delta u, \quad \mbox{in } 
\Omega\times (0,\infty)
\eq
subject to the homogeneous initial condition 
\bql{1.2}
u(\cdot\,,0)=0,\quad\mbox{in } \Omega
\eq
and the boundary conditions
\bql{1.3}
u=0,\quad \mbox{on } \Gamma_1\times[0,\infty), \quad \frac{\partial u}{\partial \nu}=g, \quad \mbox{on }\Gamma_2\times[0,\infty).
\eq
Here $\alpha$ is the thermal diffusivity,  $\nu$ denotes the outward unit normal to $\Gamma$ and $g$ is a given and sufficiently smooth function. Existence and uniqueness of classical (see \cite{Friedman,Kr}) and of
weak solutions (see \cite{La,LiMa}) of this mixed initial boundary value problem is  well established.

In this work we consider the inverse problem of determining the interior boundary curve $\Gamma_1$
from the  knowledge of the thermal flux $g$ and of  the temperature 
\bql{1.4}
u=f, \quad\mbox{on } \Gamma_2 \times (0,\infty).
\eq


\section{Laguerre transform and Boundary Integral Equations}\label{Sec2}
\setcounter{equation}{0}

We solve the system \eref{1.1} -- \eref{1.4} in two steps. First, we represent the solution $u$ as a (scaled) Fourier expansion with respect to the Laguerre polynomials resulting to a sequence of stationary mixed problems. Then, we reduce the problem to a system of boundary integral equations.

We consider the expansion
$$
 u(x,t)=\kappa\sum_{n=0}^\infty
 u_n(x)L_n(\kappa t),
 $$
 where 
$$
  u_n(x):=\int_0^\infty e^{-\kappa t}L_n(\kappa t)u(x,t)\,dt,\quad 
n=0,1,2,\ldots
$$
for $\kappa>0,$ and the Laguerre polynomial  $L_n$ of order $n,$ defined by 
$$
L_n(t) = \sum_{k=0}^n \binom{n}{k} \frac{(-t)^k}{k!},
$$
and satisfying the recurrence relations 
\begin{align*}
(n+1)L_{n+1}(t) &=(2n+1-t)L_n(t)-nL_{n-1}(t),  \\
L_{n+1}'(t) &=L_n'(t) - L_n(t), 
\end{align*}
for $n=0,1,2,\ldots.$ Using these formulas, we can show (see~\cite{ChKr2}) that the  Fourier--Laguerre coefficients $u_n$, satisfy the following sequence of mixed  problems
 \bql{lag10}
 \Delta u_n - \gamma^2 u_n
 =\beta\sum_{m=0}^{n-1} u_m ,
 \quad \mbox{ in }\Omega
 \eq
 with boundary conditions
 \bql{lag11}
 u_n=0,\quad \mbox{on }\Gamma_1\quad \textrm{and}\quad \frac{\partial u_n}{\partial \nu}= g_n ,\quad\mbox{on}\;\Gamma_2.
 \eq
 Here, 
 $$
  g_{n}(x):=\int_0^\infty e^{-\kappa t}
 L_n(\kappa t)g(x,t)\,dt,\quad n=0,1,2,\ldots
 $$
are the Laguerre--Fourier coefficients of the given function $g$ and 
 $\beta=\kappa/\alpha, \, \gamma^2 = \beta.$ The additional condition \eref{1.4} results to
\bql{lag12}
u_n=f_{n},\quad \mbox{on }\Gamma_2
\eq
with 
$$
  f_{n}(x):=\int_0^\infty e^{-\kappa t}
 L_n(\kappa t) f(x,t)\,dt,\quad n=0,1,2,\ldots.
 $$

In order to apply the non-linear integral equation method we need the fundamental solution of the equations \eref{lag10}. We consider the modified Bessel functions
 \bql{besselmod}
 I_0(z)=\sum^\infty_{n=0} \;
 \frac{1}{(n!)^2}\,\left(\frac{z}{2}\right)^{2n},
 \quad
 I_1(z)=\sum^\infty_{n=0} \;
 \frac{1}{n!(n+1)!}\,\left(\frac{z}{2}\right)^{2n+1},
 \eq
 and the modified Hankel functions
\begin{equation}\label{macdonald}
\begin{aligned}
K_0(z)  &= - \left(\ln { \frac{z}2}+C\right)\,I_0(z)
 + \sum^\infty_{n=1}
 \frac{\psi(n)}{ (n!)^2} \,\left(\frac{z}{2}\right)^{2n}, \\
 K_1(z) &= \frac{1}{z}+\left(\ln { \frac{z}2}+C\right)\,I_1(z)
 -\frac{1}{2}\sum^\infty_{n=0}
 \frac{\psi(n+1)+\psi(n)}{ n!(n+1)!}\,\left(\frac{z}{2}\right)^{2n+1}
\end{aligned} 
\end{equation}
  of order zero and one, respectively.  Here, we set $\psi(0)=0$ and
 $$
 \psi(n)=\sum_{m=1}^n\frac{1}{m}\;, \quad n=1,2,\ldots
 $$
 and $C = 0.57721\ldots$ denotes the Euler constant (see \cite{AbSt}).

 We define the polynomials $v_n$ and $w_n$ by
 $$
 v_n(r)=
 \sum_{k=0}^{\left[\frac{n}{2}\right]}a_{n,2k}r^{2k},\quad
 w_n(r)=
 \sum_{k=0}^{\left[\frac{n-1}{2}\right]}a_{n,2k+1}r^{2k+1}
 $$
 with the convention $w_0 (r) = 0.$ The coefficients are given by the relations 
\begin{align*}
a_{n,0} &=1, \\
 a_{n,n} &=-\frac{1}{2\gamma n}\;\beta a_{n-1,n-1}, \\
a_{n,k} &=\frac{1}{2\gamma k}
 \left\{4\left[\frac{k+1}{2}\right]^2a_{n,k+1}
 -\beta \sum_{m=k-1}^{n-1} a_{m,k-1}\right\},
 \quad k=n-1,\ldots,1,
\end{align*}  
  
 for $n=1,2,\ldots$. 
 
Then, the sequence of functions
 \bql{fund-sol}
 \Phi_n(x,y):=K_0(\gamma|x-y|)\,v_n(|x-y|)
 +K_1(\gamma|x-y|)\,w_n(|x-y|),\quad x\neq y,
 \eq
 satisfies  \eref{lag10} with respect to
 $x$ in $\real^2\setminus\{y\}$
 for  $n=0,1,2,\ldots,$
 that is, the function $\Phi_n$ provides a  fundamental
 solution (see \cite{ChJo}).  Now we can represent the solutions $u_n$ of the problem \eref{lag10} -- \eref{lag12} using the following single layer potential form
\bql{pot1}
 u_n(x)=\frac{1}{\pi}\sum_{\ell=1}^2\sum_{m=0}^n
 \int_{\Gamma_\ell} \phi_m^\ell(y) \Phi_{n-m}(x,y)\, ds(y),
 \quad x\in \Omega
 \eq
with the unknown densities $\phi_m^1$ and $\phi_m^2$, $m=0,1,\ldots$, defined on the boundary curves $\Gamma_1$ and $\Gamma_2$, respectively, and $\Phi_n$ is given by~\eref{fund-sol}.

We let $x$ tend to the boundary and using the boundary conditions and the standard jump relations we get the following system

\begin{alignat}{3}
\frac1{\pi} \sum_{\ell = 1}^2  \int_{\Gamma_\ell} \phi^\ell_n (y)\Phi_0 (x,y)  ds(y) &= F_{1,n} (x),   &&\quad x &&\in \Gamma_1, \label{eq1}\\
 \phi^2_n (x)+ \frac1{\pi} \sum_{\ell = 1}^2  \int_{\Gamma_\ell}\phi^\ell_n (y) \frac{\partial \Phi_0 }{\partial n (x)}  (x,y) ds(y) &= G_n (x),  &&\quad x &&\in \Gamma_2 ,\label{eq2}\\
\frac1{\pi} \sum_{\ell = 1}^2  \int_{\Gamma_\ell} \phi^\ell_n (y)\Phi_0 (x,y)  ds(y) &=  F_{2,n} (x),   &&\quad x &&\in \Gamma_2 ,\label{eq3}
\end{alignat}
with right-hand sides
\begin{align*}
F_{1,n} (x) &= -  \frac1{\pi} \sum_{\ell = 1}^2 \sum_{m=0}^{n-1} \int_{\Gamma_\ell} \phi^\ell_m (y) \Phi_{n-m} (x,y) ds(y), \\
G_n (x) &= g_n (x)  - \sum_{m=0}^{n-1} \phi^2_m (x)  -  \frac1{\pi} \sum_{\ell = 1}^2 \sum_{m=0}^{n-1} \int_{\Gamma_\ell} \phi^\ell_m (y)\frac{\partial\Phi_{n-m}}{\partial n (x)} (x,y)  ds(y), \\
F_{2,n} (x) &= f_n (x) -  \frac1{\pi} \sum_{\ell = 1}^2 \sum_{m=0}^{n-1} \int_{\Gamma_\ell} \phi^\ell_m (y)\Phi_{n-m} (x,y)  ds(y).
\end{align*}

This is a system of three  equations for the three unknowns: the two densities $\phi^1_n , \, \phi^2_n$ and the boundary $\Gamma_1.$ The operators are singular,  linear on the densities but act non-linearly on the boundary curve. We will consider the Frech\'et derivative of the integral operators for linearizing it. 

There could be many ways to solve this system of integral equations. Here, motivated by \cite{JohSle07, KrRu}, we propose the following iterative scheme which avoids the full linearization of the system and the computation of the Fr\'echet derivative of the normal derivative operator, which is more involved compared to the single layer operator. We refer to \cite{AltKre12, Greece, GinMin19} for some recent applications of this scheme to different regimes.

Before presenting the iterative method, we consider the parametrization of the system \eqref{eq1} -- \eqref{eq3}.  We assume the following parametric representation of the boundary
$$
\Gamma_\ell=\{x_\ell(s)=(x_{1\ell}(s),x_{2\ell}(s)), s\in[0,2\pi]\}, \quad \ell = 1,2
$$  
and we define
\[
\varphi_n^\ell(s) :=\phi_n^\ell(x_\ell(s))|x'_\ell(s)|.   
\]

Then, the system \eqref{eq1} -- \eqref{eq3} takes the form
\begin{alignat}{3}
\frac{1}{2\pi}\sum_{\ell=1}^2\int_0^{2\pi} \varphi_n^\ell(\sigma) H^{1,\ell}_0 (s,\sigma)\, d\sigma  &= \tilde F_{1,n} (s), &&\quad s &&\in [0,2\pi], \label{eq1par}\\
\frac{\varphi_n^2(s)}{|x'_2(s)|}+\frac{1}{2\pi}\sum_{\ell=1}^2\int_0^{2\pi} \varphi_n^\ell(\sigma) Q^{2,\ell}_0 (s,\sigma) d\sigma &= \tilde G_n (s), &&\quad s &&\in [0,2\pi], \label{eq2par}\\
\frac{1}{2\pi}\sum_{\ell=1}^2\int_0^{2\pi} \varphi_n^\ell(\sigma) H^{2,\ell}_0 (s,\sigma)\, d\sigma  &= \tilde F_{2,n} (s), &&\quad s &&\in [0,2\pi], \label{eq3par}
\end{alignat}
for $n=0,\ldots,N$, $N\in\nat$ and right-hand sides
\begin{align*}
\tilde F_{1,n} (s) &=  -\frac{1}{2\pi}\sum_{\ell=1}^2\sum_{m=0}^{n-1}\int_0^{2\pi}\varphi_m^\ell(\sigma) H^{1,\ell}_{n-m}(s,\sigma)\, d\sigma, \\
\tilde G_n (s) &= g_{n}(x_2(s))- \frac{1}{|x'_2(s)|}\sum_{m=0}^{n-1} \varphi^2_m (s)
-\frac{1}{2\pi}\sum_{\ell=1}^2\sum_{m=0}^{n-1}\int_0^{2\pi}\varphi_m^\ell(\sigma) Q^{2,\ell}_{n-m}(s,\sigma)\, d\sigma, \\
\tilde F_{2,n} (s) &= f_n (x_2 (s))  -\frac{1}{2\pi}\sum_{\ell=1}^2\sum_{m=0}^{n-1}\int_0^{2\pi}\varphi_m^\ell(\sigma) H^{2,\ell}_{n-m}(s,\sigma)\, d\sigma.
\end{align*}
The kernels are given by
\begin{equation}\label{def_par_ker}
 H^{k,\ell}_ n(s,\sigma)=2\Phi_n (x_k(s),x_\ell(\sigma)), \quad  Q^{k,\ell}_ n(s,\sigma)=2\frac{\partial \Phi_n }{\partial \nu(x)} (x_k(s), x_\ell(\sigma)),
 \end{equation} 
 for $s\neq \sigma, \,\,k,\ell=1,2$,  and $n=0,\ldots,N$. The functions $\Phi_n$ are defined in  \eref{fund-sol}.

We solve the above system of equations using the iterative scheme: 
\begin{description}
\item[Step 1] Given an initial approximation of $\Gamma_1$, we solve the sequence of well-posed systems of integral equations \eqref{eq1par} and \eqref{eq2par} for $\phi^1_n, \, \phi^2_n, \,\, n = 0,...,N.$
\item[Step 2] Keeping now the densities fixed,  we linearize  the ill-posed integral equation \eqref{eq3par} resulting to
\begin{equation}\label{eq_frechet}
 \sum_{m=0}^n  \mathcal D_{n-m} [\varphi_m^1, x_1 ; \chi] (s) = f_n (x_2 (s)) - \frac{1}{2\pi} \sum_{\ell = 1}^2 \sum_{m=0}^n \int_0^{2\pi} \varphi_m^\ell (\sigma) H^{2,\ell}_{n-m} (s,\sigma)\, d\sigma.
\end{equation}
We solve the $n = 0,...,N$  equations for the parametrization $\chi$ of the perturbed  $\Gamma_1,$ and we update as $x_1 + \chi .$
\end{description}

Equation \eqref{eq_frechet} contains the Fr\'echet derivative $\mathcal D_n$ of the integral operator with kernel $H^{2,1}_n$ with respect to $x_1$ as a linear operator on $\chi.$ The analytic form of this operator is given in the next section. In \eqref{eq_frechet}  we have $N$ equations to be solved for one function $\chi.$ We present later two different ways for solving this equation.

\section{Numerical Implementation}\label{sec3}
\setcounter{equation}{0}

In this section we consider the numerical implementation of the iterative scheme presented above. We refer to \eqref{eq1par} -- \eqref{eq2par} as the ``field" system and to \eqref{eq3par} as the ``data" equation.

\subsection{Numerical solution of the ``field" system}\label{sec3.1}

The singular kernels $H^{\ell,\ell}_n$ appearing in \eqref{eq1par} and in \eqref{eq3par}, are analyzed using the explicit expression \eref{fund-sol} for the elements of the fundamental sequence together with the expansions \eref{besselmod} --  \eref{macdonald}. Then, we get the expression
\begin{align*}
 H^{k,\ell}_n(s, \sigma) 
& = K_0(\gamma|x_k (s)-x_\ell (\sigma)|)v_n(|x_k (s)-x_\ell (\sigma)|) \\
&\phantom{=}+K_1(\gamma|x_k (s)-x_\ell (\sigma)|)w_n(|x_k (s)-x_\ell (\sigma)|).
 \end{align*}

After lengthy but straightforward calculations, we treat the logarithmic singularity, using the following representation
$$
 H^{\ell,\ell}_n(s,\sigma)= H^{\ell,\ell}_{n,1}(s,\sigma)\ln \left( \frac{4}{e}\sin^2 \frac{s-\sigma}{2}\right)+H^{\ell,\ell}_{n,2}(s,\sigma),
$$
 where
 \begin{align*}
 H^{\ell,\ell}_{n,1}(s, \sigma) 
& = -I_0(\gamma|x_\ell (s)-x_\ell (\sigma)|)v_n(|x_\ell (s)-x_\ell (\sigma)|) \\
&\phantom{=}+I_1(\gamma|x_\ell (s)-x_\ell (\sigma)|)w_n(|x_\ell (s)-x_\ell (\sigma)|) 
 \end{align*}
 and
$$
 H^{\ell,\ell}_{n,2}(s,\sigma)=H^{\ell,\ell}_n(s,\sigma) - H^{\ell,\ell}_{n,1}(s,\sigma)\ln \left( \frac{4}{e}\sin^2 \frac{s-\sigma}{2}\right).
$$
The diagonal terms are given by
$$
H^{\ell,\ell}_{n,2}(s, s)=-2C-1-2\ln\left(\frac{\gamma|x'_\ell(s)|}{2}\right)+\frac{2a_{n,1}}{\gamma},\quad n=0,1,\ldots,N.
$$
 
For the representation of the kernels $Q^{k,\ell}_n$, we introduce the function
 $$
 h_{k,\ell}(s,\sigma)=\frac{(x_{k,1}(s)-x_{\ell,1} (\sigma))x'_{k,2}(s)-(x_{k,2}(s)-x_{\ell,2} (\sigma))x'_{k,1}(s)}{|x_k(s)-x_\ell (\sigma)|}
 $$
and the polynomials
\begin{align*}
\tilde{v}_n(r) &=\gamma \sum_{m=0}^{\left[\frac{n}{2}\right]}a_{n,2m}r^{2m}-2 \sum_{m=1}^{\left[\frac{n-1}{2}\right]}ma_{n,2m+1}r^{2m}, \\
 \tilde{w}_n(r) &=\gamma \sum_{m=0}^{\left[\frac{n-1}{2}\right]}a_{n,2m+1}r^{2m+1}-2 \sum_{m=1}^{\left[\frac{n}{2}\right]}ma_{n,2m}r^{2m-1}.
\end{align*}

 Then, the kernels $Q^{k,\ell}_n$ admit the form
\begin{align*}
 Q^{k,\ell}_n (s,\sigma) &=  
2 h_{k,\ell}(s,\sigma)\,\{K_1(\gamma|x_k(s)-x_\ell (\sigma)|)
 \tilde v_n(|x_k (t)-x_\ell (\sigma)|) \\
&\phantom{=}+
 K_0(\gamma|x_k(s)-x_\ell (\sigma)|)
 \tilde w_n(|x_k(s)-x_\ell (\sigma)|)\}
\end{align*} 
 for $s\neq\sigma$ and $n=0,1,\ldots,N$. 

The kernels $Q^{\ell,\ell}_n$ have logarithmic singularity. As in the case of the kernels $H^{\ell,\ell}_n$ performing similar calculations,  we derive the decomposition
$$
 Q^{\ell,\ell}_n(s,\sigma)= Q^{\ell,\ell}_{n,1}(s,\sigma)\ln \left( \frac{4}{e}\sin^2 \frac{s-\sigma}{2}\right)+Q^{\ell,\ell}_{n,2}(s,\sigma),
$$
 where
 \begin{align*}
 Q^{\ell,\ell}_{n,1}(s, \sigma) &=  
2 h_{\ell,\ell}(s,\sigma)\,\{I_1(\gamma|x_k(s)-x_\ell (\sigma)|)
 \tilde v_n(|x_k (s)-x_\ell (\sigma)|)\\ 
 &\phantom{=}-
 I_0(\gamma|x_k(s)-x_\ell (\sigma)|)
 \tilde w_n(|x_k(s)-x_\ell (\sigma)|)\}
\end{align*}  
 and
$$
 Q^{\ell,\ell}_{n,2}(s,\sigma)=Q^{\ell,\ell}_n(s,\sigma) - Q^{\ell,\ell}_{n,1}(s,\sigma)\ln \left( \frac{4}{e}\sin^2 \frac{s-\sigma}{2}\right).
$$
The diagonal elements are now given by
$$
 Q^{\ell,\ell}_{n,2}(s, s)=\frac{x^\prime_{\ell,2}(s)x^{\prime\prime}_{\ell,1}(s)-x^\prime_{\ell,1}(s)x^{\prime\prime}_{\ell,2}(s)}
 {|x_\ell^\prime(s)|^3},\quad n=0,1,\ldots,N.
$$

Now we have the explicit representations of the singular kernels, and we can apply the following standard quadrature rules \cite{Kr} for the numerical discretization 
\begin{align*}
 \frac{1}{2\pi}
 \int_0^{2\pi}
 f(\sigma)\,
 d\sigma
 &\approx 
 \frac{1}{2M}\sum_{k=0}^{2M-1}
  f(s_k), \\
   \frac{1}{2\pi}
 \int_0^{2\pi}
 f(\sigma)
 \ln\left(\frac{4}{e} \sin^2 \frac{s- \sigma}2 \right)
 d\sigma
 &\approx
 \sum_{k=0}^{2M-1}
 {R}_{k}(s)\,f(s_k),
\end{align*}
for the mesh points 
\begin{equation}\label{mesh_points}
s_k=kh, \quad k=0,\ldots,2M-1, \quad h=\pi/M, \quad M \in \nat,
\end{equation} 
using the weight function 
$$
 {R}_{k}(s) = 
\displaystyle
 - \frac{1}{2M} \;  \left(1+2\sum^{M-1}_{m=1} \frac{1}{ m} \,
\cos m(s- s_k) - \frac1{ M} \, \cos M(s-s_k )\right).
$$

Thus, we approximate the solution of the ``field" system \eqref{eq1par} -- \eqref{eq2par} by collocating the integral equations at the nodal points $\{s_k\}$ leading to the sequence of linear systems
$$
\begin{aligned}
\sum\limits_{j = 0}^{2M-1}\left\{\varphi_{n,j}^1\left[R_j(s_i){H}^{1,1}_{0,1}(s_i, s_j)+\frac{1}{2M}{H}^{1,1}_{0,2}(s_i, s_j)\right]+
\varphi_{n,j}^2\frac{1}{2M} H^{1,2}_0 (s_i, s_j) \right\} &= \tilde F_{1,n} (s_i), \\
\sum\limits_{j = 0}^{2M-1} \left\{\varphi_{n,j}^1 \frac{1}{2M}Q^{2,1}_0(s_i, s_j) + 
	\varphi_{n,j}^2\left[R_j(s_i)Q^{2,2}_{0,1}(s_i, s_j)+\frac{1}{2M} Q^{2,2}_{0,2}(s_i, s_j)\right]\right\} \\
	+ \frac{\varphi_{n,i}^2}{|x_2'(s_i)|} &= \tilde G_n (s_i),
\end{aligned}
$$
 for $i=0,\ldots,2M-1$, with the right-hand sides
\begin{align*}
\tilde F_{1,n} (s_i) &= 
 -\sum_{j=0}^{2M-1}\sum_{m=0}^{n-1}\left\{\varphi_{m,j}^1 [R_j(s_i){H}^{1,1}_{n-m,1}(s_i, s_j)+\frac{1}{2M}{H}^{1,1}_{n-m,2}(s_i, s_j)] \right.\\
 &\phantom{=} +
  \left.\varphi_{m,j}^2 \frac{1}{2M}H^{1,2}_{n-m}(s_i, s_j)\right\}\nonumber
\end{align*}
 and
 \begin{align*}
\tilde{G}_n (s_i) &= g_{n}(x_2(s_i))-\frac{1}{|x_2'(s_i)|}\sum_{m=0}^{n-1}\varphi_{m,i}^2
 -\sum_{j=0}^{2M-1}\sum_{m=0}^{n-1}\left\{\varphi_{m,j}^1 \frac{1}{2M}Q^{2,1}_{n-m}(s_i, s_j)\right.\\
 &\phantom{=} +
  \left. \varphi_{m,j}^2 [R_j(s_i) Q^{2,2}_{n-m,1}(s_i, s_j)+\frac{1}{2M}Q^{2,2}_{n-m,2}(s_i, s_j)]\right\},
 \end{align*}
 where we used the abbreviation $\varphi_{n,j}^\ell \approx \varphi_n^\ell(s_j)$, $\ell = 1,2$, $n=0,\ldots,N$, $j=0,\ldots,2M-1$.

\subsection{Numerical solution of the ``data" equation}\label{sec3.2}

As described in our algorithm, we search for the correction of $\Gamma_1$ by solving the ``data" equation \eqref{eq3par}, assuming that we know the densities 
$\varphi^\ell_n$, $\ell=1,2$, $n=0,...,N$.

For simplicity,  we consider starlike interior curve, meaning we assume  parametrization in polar coordinates of the form
$$
x_1(s)=\{r(s) (\cos s,\sin s) : \, s\in  [0,2\pi]\},
$$
where $r:\real \to(0,\infty)$ is 
a  $2\pi$ periodic function representing the radial distance from the origin. However, we wish to stress that the following analysis is also applicable to other boundaries.

The linearized equation \eqref{eq_frechet} admits the following parametric form
\begin{equation}\label{lin_eq}
\sum_{m=0}^{n} \mathcal D_{n-m} [\varphi^1_m, r;q]  (s)  =  f_n (x_2 ( s)) -  \frac1{2\pi} \sum_{\ell=1}^2 \sum_{m=0}^{n}\int_0^{2\pi}  \varphi^\ell_m (\sigma) H^{2,\ell}_{n-m} (s,\sigma)  d\sigma, 
\end{equation}
where $q$ is the radial function of the perturbed boundary. The Fr\'echet derivative $\mathcal{D}_n$ has the explicit form
\[
 \mathcal D_n [\varphi , r;q]  (s) = \frac1{2\pi} \int_0^{2\pi} q (\sigma)   \varphi (\sigma) D_n (s,\sigma) d\sigma
\]
with kernel
$$
 D_n(s,\sigma) = - \frac{(x_2 (s) -  x_1 (\sigma) ) \cdot (\cos \sigma,\,\sin \sigma)}{|x_2 (s) -  x_1 (\sigma) |}   \tilde{\Phi}_n (|x_2 (s) -  x_1 (\sigma) |),
$$
where
$$
 \tilde{\Phi}_n (r) = K_1 (\gamma r)\,\tilde{v}_n (r)  +K_0 (\gamma r)\,\tilde{w}_n ( r).
$$

\begin{theorem}
The Fr\'echet  derivative operator $\mathcal D_n [\varphi , r;q] $ is injective at the exact solution.
\end{theorem}

\proof
We just have to show that if $q$ solves
\[
\sum_{m=0}^{n} \mathcal D_{n-m} [\varphi^1_m, r;q]  (s)  =  0,
\]
then $q=0.$
We follow the ideas of \cite{ChIv} and we set
$$
V_n(x)=\frac{1}{\pi}\sum_{m=0}^{n} \int_{\Gamma_1}\varphi_m^1(y)(\zeta(y),\partial_y)\Phi_{n-m}(x,y)ds(y)\quad x\in\real^2\setminus \Gamma_1,
$$
where $\zeta(x_1(s))=q(s)(\cos s,\sin s)$. Clearly $V_n$ satisfies the sequence \eref{lag10}. We define $V_n^{\pm}(x) := \lim_{h\rightarrow 0^+} V_n (x \pm h \nu).$ Then, by 
assumption $V_n^+|_{\Gamma_1}=0$.
The perturbed interior curve can be represented as follows \cite{ChIv,IvKr}
$$
\Gamma_{1,r+q}=\{r(s)(\cos s,\sin s)+\tilde{q}(s)\nu(x_1(s)): s\in[0,2\pi] \},
$$
for small perturbations and a given function $\tilde{q}.$

Then we can rewrite the functions $V_n$ in the form
$$
V_n(x)=\frac{1}{\pi}\sum_{m=0}^{n} \int_{0}^{2\pi}\tilde{q}(\sigma) \varphi _m^1(\sigma)(\nu(x_1(\sigma)),\partial_{x_1(\sigma)})\Phi_{n-m}(x,x_1(\sigma))|x_1'(\sigma)|d\sigma, \quad x\in\real^2\setminus \Gamma_1.
$$
By the properties of the double layer potentials $V_n$ \cite{ChKr3},  we extend it continuously to $\Gamma_1$ as 
\begin{align*}
V_n^\pm (x_1(s)) &=\pm \tilde{q}(s) \sum_{m=0}^{n} \varphi_m^1(s) \\
&\phantom{=}+\frac{1}{\pi}\sum_{m=0}^{n} \int_{0}^{2\pi}\tilde{q}(\sigma)\varphi_m^1(\sigma)(\nu(x_1(\sigma)),\partial_{x_1(\sigma)})\Phi_{n-m}(x_1(s),x_1(\sigma))|x_1'(\sigma)|d\sigma.
\end{align*}
By the uniqueness of the exterior and interior Dirichlet problems we have
\begin{equation}\label{eq_phi}
\tilde{q}(s) \sum_{m=0}^{n} \varphi_m^1(s)=0, \quad s\in[0,2\pi].
\end{equation}

The functions $u_n$ given by \eref{pot1} solve the Dirichlet problem in the interior of $\Gamma_1$ with homogeneous boundary conditions. Then by the  unique solvability, the functions $u_n$ have to vanish in the interior of $\Gamma_1$ and hence $\tfrac{\partial u_n^-}{\partial \nu}=0$ on $\Gamma_1$. The jump relations imply that $\frac{\partial u_n^+}{\partial \nu}|_{\Gamma_1}=\sum_{m=0}^{n}  \varphi_m^1(s)$. Employing Holmgren's uniqueness theorem to the sequence of Helmholtz equations, one can show that the
Cauchy data $(u_n^+, \, \frac{\partial u_n^+}{\partial \nu})$ cannot be identically zero on an open subset of $\Gamma_1$. Thus, in view of \eqref{eq_phi}, we conclude that $\tilde{q} =0$ and therefore $q=0.$
\beweisende

We apply the quadrature rules, analyzed in the previous section, in \eqref{lin_eq} and then a collocation method approximating the function $q$ by a trigonometric polynomial of the form 
\begin{equation}\label{eq_radial}
q(s) \approx \sum_{j=0}^{2J} q_j \tau_j (s), \quad \nat \ni J \ll M,
\end{equation}
with
\[
\tau_j (s) =  \left.
  \begin{cases}
    \cos (js), & \text{for } j = 0,...,J ,\\
    \sin ((j-J)s), & \text{for } j = J+1,...,2J,
  \end{cases}
  \right.
\]

We substitute \eqref{eq_radial} in \eqref{lin_eq} and at the nodal points $\{s_i\}$ we obtain the following linear system 
\begin{equation}\label{eq_linear}
\fet A (n) \fet q = \fet b (n), \quad n = 0,...,N,
\end{equation}
where $\fet q  = (q_0, ..., q_{2J})^\top \in \real^{2J+1},$ and $\fet A \in \real^{(2M)\times (2J+1)}, \, \fet b \in \real^{2M}$ given by
\begin{align*}
\fet A_{ij} (n) &= \frac{1}{2M}\sum_{k=0}^{2M-1} \tau_j(s_k)\sum_{m=0}^{n}  \varphi^1_{m,k} L_{n-m}(s_i, s_k), \\
\fet b_i (n) &= f_n (x_2(s_i)) -\sum_{k=0}^{2M-1}\sum_{m=0}^n\left\{ \varphi_{m,k}^1 \frac{1}{2M}H^{2,1}_{n-m}(s_i,s_k) \right. \\
&\left.\phantom{=}+ \varphi_{m,k}^2\left[R_k (s_i) H^{2,2}_{n-m,1}(s_i,s_k)+\frac{1}{2M}H^{2,2}_{n-m,2}(s_i,s_k)\right]\right\}.
\end{align*}
We solve the linear system \eqref{eq_linear} either for $n=N,$ or using all available information, meaning
\[
\begin{pmatrix}
\fet A (0) \\
\fet A (1) \\
\vdots \\
\fet A (N)
\end{pmatrix} \fet q = 
\begin{pmatrix}
\fet b (0) \\
\fet b (1) \\
\vdots \\
\fet b (N)
\end{pmatrix}.
\]
We refer to these cases as the ``final-step" and the ``multi-step" system, respectively.  

Due to the ill-posedness of \eref{eq_linear} and its overdetermination we apply the least-squares method with Tikhonov regularization. Then, we solve
\[
\min_{\fet q} \{  \| \fet A \fet q - \fet b\|_2^2 + \lambda  \|  \fet q \|_2^2 \},
\] 
with the regularization parameter $\lambda >0,$ to be chosen by trail and error.

\section{Numerical Results}\label{sec4}
\setcounter{equation}{0}

We present numerical examples for different boundary curves. We consider the following three cases:
\begin{description}
\item[Example 1] The exterior boundary curve $\Gamma_2$ is a circle with center $(0,0)$ and radius $1,$ and the interior boundary curve $\Gamma_1$ (to be reconstructed) is peanut-shaped with radial function
\[
r(s) =  \sqrt{(0.5 \cos s)^2 + (0.25 \sin s)^2}.
\]
\item[Example 2] The exterior boundary curve $\Gamma_2$ is a rounded rectangle with radial function
\[
r_2 (s) = (\cos^{10} s  +\sin^{10} s )^{-0.1}
\]
and $\Gamma_1$ is a apple-shaped boundary with radial function
\[
r_1 (s) =  \frac{0.45 + 0.3 \cos s - 0.1 \sin 2s}{1+0.7 \cos s}.
\]
\item[Example 3]  Both boundary curves are kite-shaped with parametrizations
\[
x_1 (s) = \frac13 ( \cos s + 0.55 \cos 2s - 0.5, \, 1.2 \sin s)
\]
and
\[
x_2 (s) =  ( \cos s + 0.8 \cos 2s - 0.5, \, 1.5 \sin s).
\]
\end{description}

We generate the simulated Cauchy data by solving the sequence  
\[
 \Delta u_n - \gamma^2 u_n
 =\beta\sum_{m=0}^{n-1} u_m,
 \quad \mbox{ in }D
\]
 with boundary conditions
 \[
 u_n= f_{1,n}, \quad \mbox{on }\Gamma_1,\quad \text{and}\quad  u_n = f_{2,n},\quad\mbox{on}\;\Gamma_2,
 \]
for given boundary functions $f_{\ell,n}, \, \ell = 1,2.$ To avoid an inverse crime, we consider double amount of nodal points compared to the inverse problem and we compute the thermal flux $g$ by tending the solution to the exterior boundary and considering the jump relations. Then, we add error to the Cauchy data with respect to the $L^2$ norm
\[
f^\delta_{2,n} = f_{2,n} + \delta \frac{\|f_{2,n} \|_2}{\|u\|_2} u,  \quad \text{and} \quad g^\delta_{2,n} = g_{2,n} + \delta \frac{\|g_{2,n} \|_2}{\|v\|_2} v,
\]
for given noise level $\delta,$ and normally distributed random variables $u,v \in \real.$ At every step, we update the regularization parameter using
\[
\lambda_k = \lambda_0 \,0.9^{k-1}, \quad k = 1,2,... 
\]

The boundary functions are chosen as
\[
f_{1,n}=0, \quad \text{and} \quad f_{2,n} = \frac{e (2+\kappa n (\kappa (n-1)-4))}{4 (\kappa +1)^{n+3}}, \quad n = 0,...,N.
\]

\begin{figure}[t!]
  \centering
  \includegraphics[scale=0.8]{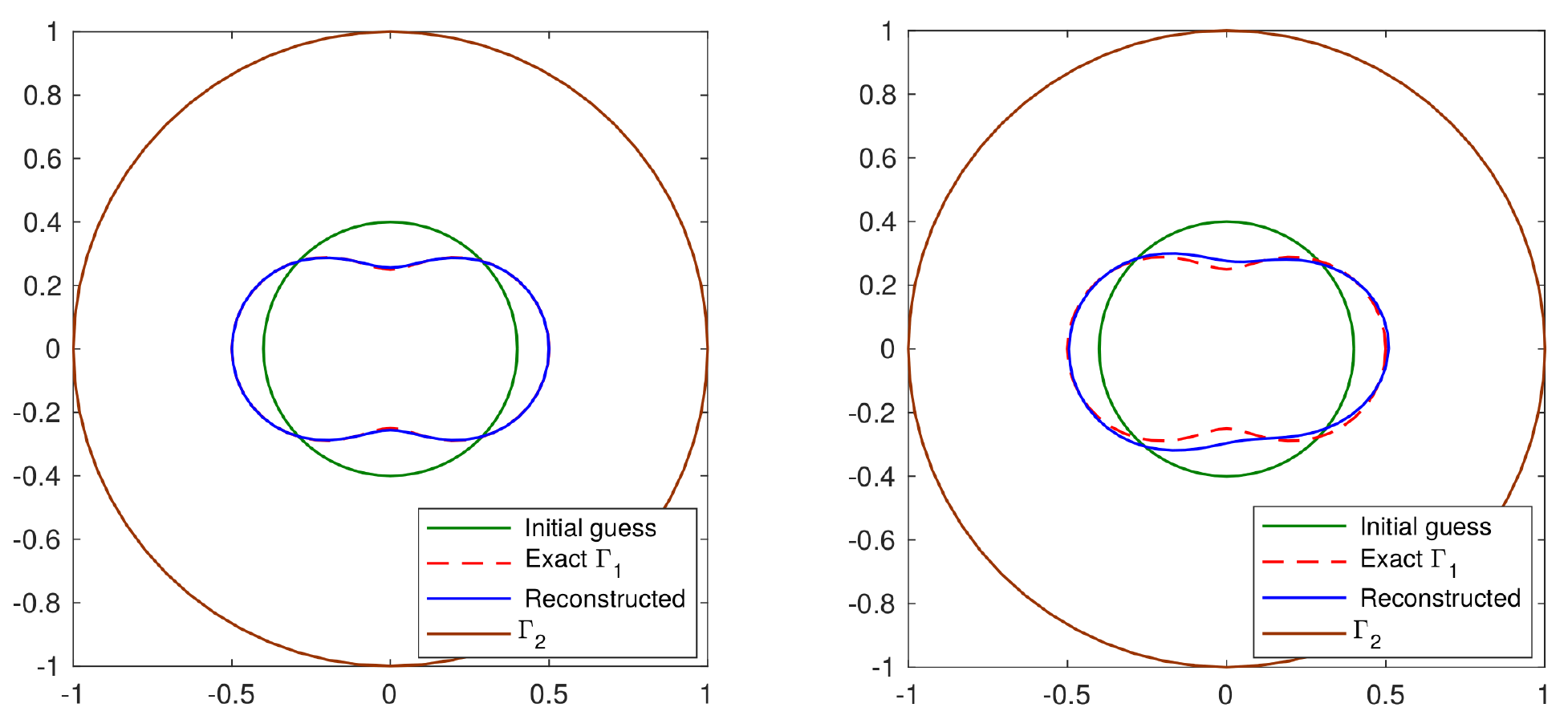}
\caption{Reconstructions of the peanut-shaped boundary $\Gamma_1$ for exact data (left) and data with $3\%$ noise (right). }\label{Fig1}
\end{figure}

In all examples, the initial guess is a circle with center $(0,0)$ and radius $r_0 ,$ we set $\alpha =1$ and  $M=64.$ We observed that the additional information provided by the ``multi-step" system produces more stable reconstructions for noisy data and requires less iterations for exact data. Thus, in the following examples we consider the ``final-step" method for exact data and the ``multi-step" method for noisy data.
In the first example, we set $\kappa =1$ and $N=10. $ The results are presented in \autoref{Fig1}, for $J=5$ coefficients, using $r_0 = 0.4$ as initial radius. The reconstruction for noise-free data are obtained using as initial regularization parameter $\lambda_0 = 0.0001$  and 24 iterations. For data with $3\%$ noise, we used $\lambda_0 = 0.001$ and 14 iterations. 

\begin{figure}[t!]
  \centering
  \includegraphics[scale=0.75]{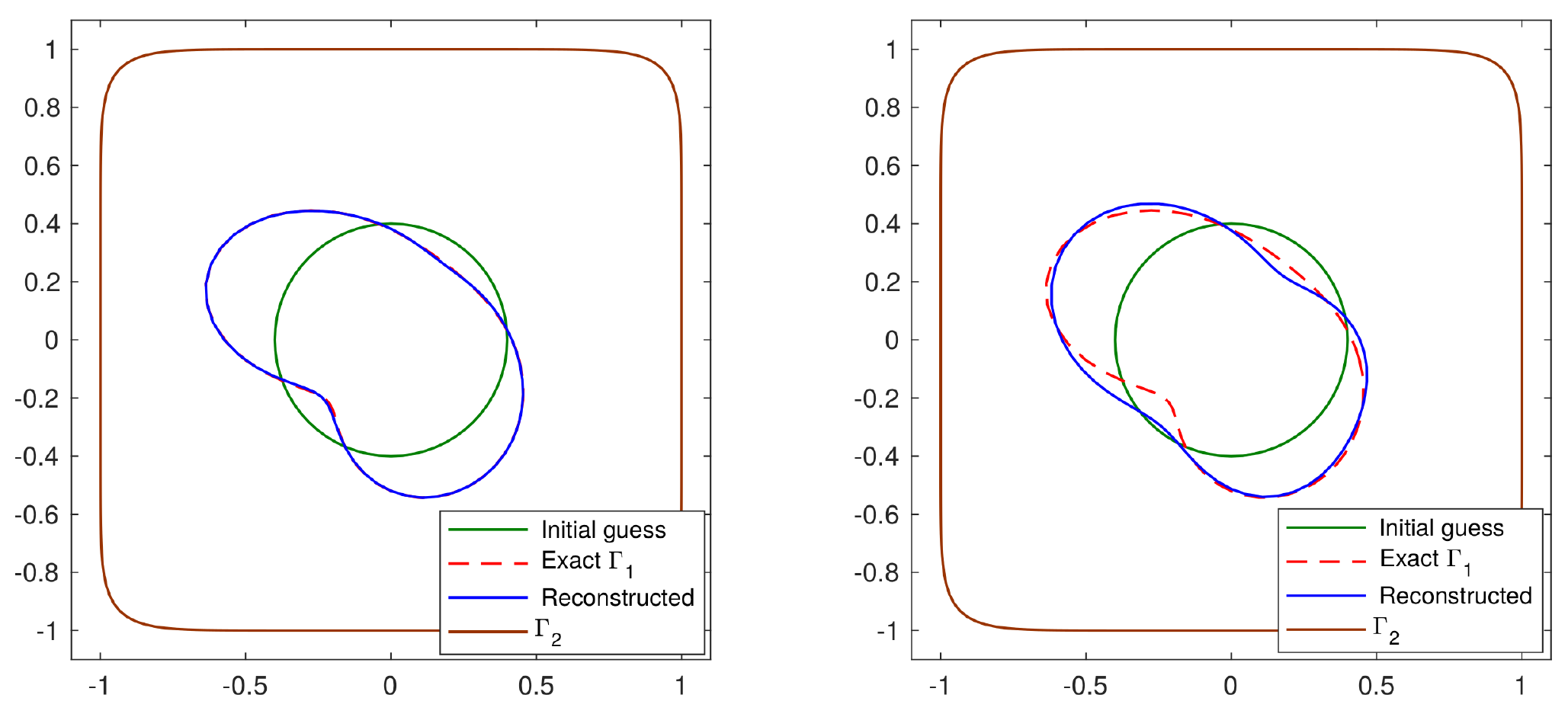}
\caption{Reconstructions of the apple-shaped boundary $\Gamma_1$ for exact data (left) and noisy data (right). }\label{Fig2}
\end{figure}

In \autoref{Fig2}, we see the reconstructions for the parametrized boundary curves of Example 2. We keep all the parameters the same as in the first example. We mark here that the results are again satisfactory independently of the boundary parametrizations. The regularization parameters are also kept fixed and the reconstructions are after 20 iterations (exact data) and after 9 iterations (data with $3\%$ noise). 

\begin{figure}[t!]
  \centering
  \includegraphics[scale=0.78]{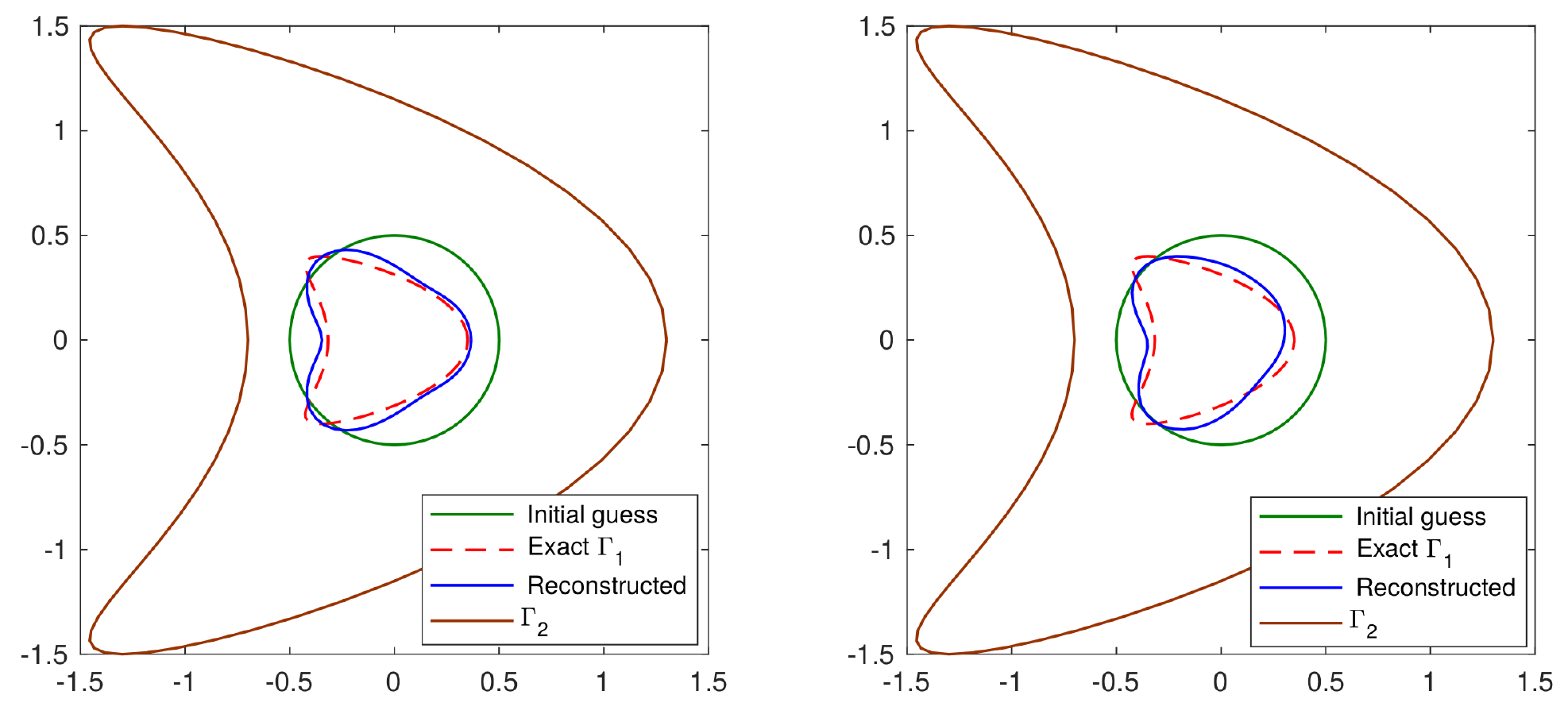}
\caption{Reconstructions of the kite-shaped boundary $\Gamma_1$ for exact data (left) and data with $3\%$ noise (right). }\label{Fig3}
\end{figure}

In the third example, we test the performance of the algorithm for more irregular boundary curves. Here, we set $J=7$ and the initial radius is $r_0 = 0.5.$ The results presented in \autoref{Fig3} are after 13 iterations for exact data and after 9 iterations for noisy data. 

\section{Conclusion}
We developed a non-linear integral equations approach for the inverse parabolic problem related to the reconstruction of a part of the boundary curve from the knowledge of the Cauchy data on the other part of the  boundary. Our strategy consisted on the consecutive dimension reduction. Firstly we considered the semi-discretization by Laguerre transform in time for the given three-dimensional problem. 
It resulted to a sequence of stationary inverse boundary problems for the Helmholtz equation. With the help of the modified single layer potentials representation, these problems were reduced to a sequence of one-dimensional non-linear boundary integral equations. 
Then, a Newton-type iteration method was applied. We solved the well-posed system of linear integral equations by the Nystr\"om method and the ill-posed linear integral equation by the collocation method with Tikhonov regularization, at every iteration step. In general, our approach can be applied without significant changes to the inverse boundary problems for the hyperbolic equation and can also be extended to the case of three-dimensional domains.

\end{document}